\documentclass[12pt, a4paper]{amsart}

\usepackage[utf8]{inputenc}
\usepackage[T1]{fontenc}
\usepackage[english]{babel}
\usepackage[a4paper, margin=2.5cm]{geometry}
\usepackage{amssymb,mathtools,mathrsfs}
\usepackage{enumitem}
\usepackage{xcolor}

\usepackage{hyperref}
\hypersetup{
    colorlinks=true,
    linkcolor=blue,
    citecolor=red,
    urlcolor=blue
}
\usepackage{cleveref}

\theoremstyle{plain}
\newtheorem{theorem}{Theorem}[section]
\newtheorem{lemma}[theorem]{Lemma}
\newtheorem{corollary}[theorem]{Corollary}
\newtheorem{proposition}[theorem]{Proposition}

\theoremstyle{definition}

\newtheorem{example}[theorem]{Example}
\newtheorem{counter-example}[theorem]{Counter-example}

\theoremstyle{remark}
\newtheorem{remark}[theorem]{Remark}

\newcommand{\R}{\mathbb{R}}

\newcommand{\Mbar}{\overline{M}}

\newcommand{\Ricbar}{\overline{\operatorname{Ric}}}

\newcommand{\Hess}{\operatorname{Hess}}

\numberwithin{equation}{section}
\newcommand{\nbar}{\overline\nabla}
\newcommand{\Ric}{\operatorname{Ric}}

\begin{document}

\title[Rigidity of spacelike graphs in generalized Robertson--Walker spacetimes]{Rigidity of spacelike graphs in generalized Robertson--Walker spacetimes via the infinity-Laplacian}

\author[E. A. Lima Jr., G. Nascimento, A. Romero]{E. A. Lima Jr.$^{1}$, G. Nascimento$^{2}$, A. Romero$^{3}$}

\address{$^{1,2}$Departamento de Matem\'atica, Universidade Federal da Para\'iba, Jo\~ao Pessoa 58.051-900, Para\'iba, Brazil.}
\email{eraldo.lima@academico.ufpb.br}
\email{razielgfn@gmail.com}

\address{$^3$Departamento de Geometr\'ia y Topolog\'ia, Facultad de Ciencias,\newline Universidad de Granada, Avenida Fuente Nueva S.N., 18071 Granada, Spain.}
\email{aromero@ugr.es}


\begin{abstract} In the search for new characterizations of spacelike slices in generalized Robertson–Walker spacetimes, we introduce a new approach. The main tool for addressing this question is the $\infty$-Laplacian with respect to the Riemannian metric induced on the spacelike hypersurface, together with a Liouville-type property for the $\infty$-Laplacian on any complete Riemannian manifold \cite{AraujoMariPessoa2021}. Starting from the $\infty$-Laplacian of the height function of a spacelike hypersurface, we derive several differential equalities involving the $\infty$-Laplacian of suitable compositions of the height function, namely, the warping function $f(h)$, its primitive $P(h)=\int^h f$, and the logarithmic composition $\log(f(h))$. Under some natural assumptions, including timelike and null energy conditions, global rigidity results are established in the parametric setting. Some of these results specialize naturally to complete spacelike hypersurfaces in the Steady State spacetime and in Einstein–de Sitter spacetime. Finally, we consider a second-order PDE involving the $\infty$-Laplacian and study the uniqueness of its entire solutions on the fiber of the generalized Robertson-Walker spacetime. This leads to new Calabi–Bernstein-type uniqueness results in the non-parametric setting.
\end{abstract}

\maketitle
\footnote[1]{The first author was partially supported by CNPq, Brazil, through PQ-C Grant 305876/2025-2 and Universal Grant 406078/2025, and by FAPESQ, Brazil, through its Universal program (Grant Agreement No.~2453/2026)}
\footnotetext[3]{The third author was partially supported by the Spanish MICINN and ERDF project PID2024-156031NB-I00.}

\thispagestyle{empty}

\section{Introduction}\label{sec:introduction}

\hyphenation{mo-dels}
\hyphenation{co-mo-ving}
\hyphenation{geo-me-tric}
\hyphenation{cons-tant}
\hyphenation{ri-gi-di-ty}

Since the early development of General Relativity, Robertson-Walker (or Friedmann-Lemaître-Robertson-Walker) spacetimes have been considered suitable mathematical models to describe an isotropic and spatially homogeneous relativistic universe on a large scale. However, on a smaller scale, these properties are less realistic \cite{RS1995}. Thus, the notion of Generalized Robertson-Walker (GRW) spacetime was introduced \cite{AliasRomeroSanchez1995} and it has been shown to be useful not only from a mathematical viewpoint.

\vspace{1mm}

Each GRW spacetime is foliated by a distinguished family of spacelike hypersurfaces: its spacelike slices. Each slice can be regarded as the physical space relative to a comoving observer at a fixed instant of its proper time. Under suitable geometric assumptions, with realistic physical interpretations, one can conclude that a certain complete spacelike hypersurface must be a spacelike slice. One of the  most important assumptions is the constancy of the mean curvature of the spacelike hypersurface \cite{ChoquetBruhat1976}. This rigidity result carries an analytical counterpart: the so-called Calabi-Bernstein problem. Namely, the determination of all entire spacelike graphs of constant mean curvature, is a question that lies in the realm of PDEs on Riemannian manifolds.

\vspace{1mm}

Although this paper is devoted to studying the global geometry of spacelike hypersurfaces in GRW spacetimes, our approach is partly inspired by the central role played by trapped surfaces in General Relativity and Lorentzian geometry. In fact, since Penrose's pioneering work on gravitational collapse and spacetime singularities \cite{Penrose1965}, trapped surfaces and their higher-dimensional generalizations have attracted considerable research interest (see, for instance, \cite{Senovilla2011, AliasCanovasColares2017, AndradeCruzFigueiraLima2025}). The geometric content of trappedness is encoded in the causal character of the mean curvature vector field and, more generally, in the causal orientation of the second fundamental form. This viewpoint naturally leads to rigidity and non-existence problems: under suitable assumptions on energy conditions, completeness, and boundedness, one may ask whether a complete spacelike submanifold can be genuinely tilted in time or whether it must be contained in a spacelike slice \cite{AndradeCruzFigueiraLima2025}.  

\vspace{1mm}

Our main tool for addressing this question is the infinity-Laplacian computed with respect to the Riemannian metric induced on the spacelike hypersurface, namely,
\begin{equation*}
  \Delta_\infty u=\Hess u(\nabla u,\nabla u).
\end{equation*}
The key analytic input is the Liouville-type property for the infinity-Laplacian on any complete Riemannian manifold, obtained from the potential-theoretic characterization of completeness by Araujo, Mari and Pessoa \cite{AraujoMariPessoa2021}. When applied to the height function of a spacelike hypersurface in a GRW spacetime, this principle turns pointwise geometric information into global slice rigidity. Thus, our starting point is closely related to the recent work of Andrade, Cruz, Figueira and Lima \cite{AndradeCruzFigueiraLima2025}, where rigidity and non-existence results for complete totally trapped spacelike submanifolds in GRW spacetimes are obtained by means of the aforementioned Liouville property for the infinity-Laplacian. 

\vspace{1mm}

Among our parametric rigidity results is Corollary \ref{cor:directional-zero}, which states:
\begin{quote}{\it 
Let $\psi: M^n\to \Mbar$ be a complete spacelike hypersurface in a GRW spacetime $\Mbar$, and let $h=\pi_I\circ\psi$ be the height function on $M$. Suppose that the warping function $f$ of $\Mbar$ and the Weingarten operator $A$ of $\psi$ satisfy 
\begin{equation*}
  (\log f)''(h)\leq0
  \qquad\text{and}\qquad
  \langle A\nabla h,\nabla h\rangle=0
\end{equation*}
on $M$. Then the composition $f(h)$ is constant. Moreover, if $f$ is not locally constant on $I$, then $h\equiv t_0$ for some $t_0\in I$, and we have $\psi(M)=\{t_0\}\times F$.}
\end{quote} 

The convexity of the function $-\log f$ is a quite natural assumption and is related to certain energy conditions imposed on $\Mbar$. In fact, it is automatically satisfied if $\Mbar$ satisfies the Timelike Convergence Condition (TCC). On the other hand, if $\Mbar$ satisfies the Null Convergence Condition (NCC) and the Ricci tensor of the fiber $F$ of $\Mbar$ satisfies $\Ric^F\leq 0$, then \eqref{eq:ncc-grw} implies that $(\log f)''\leq 0$.
Furthermore, hypothesis $\langle A\nabla h,\nabla h\rangle=0$ is directional. Notice that this condition clearly holds for any spacelike slice. On the other hand, since we are not making any assumption on the zero set of $\nabla h$, this second hypothesis is no longer related to asymptotic directions on the spacelike hypersurface. 

\vspace{1mm}

We also obtain the following obstruction result in Corollary \ref{cor:tcc-log-rigidity}:

\begin{quote}{\it 
Let $\psi: M^n\to \Mbar$ be a complete spacelike hypersurface in a GRW spacetime $\Mbar$. Assume $\Mbar$ satisfies the TCC and 
\begin{equation*}
  f(h)\geq c>0,
\end{equation*}
for some constant $c>0$. Suppose that one of the following alternatives holds:
\begin{enumerate}
  \item[\rm(1)] $\psi$ is totally future trapped and $f'(h)\geq0$ on $M$,

\vspace*{1mm}
  
  \item[\rm(2)] $\psi$ is totally past trapped and $f'(h)\leq 0$ on $M$.
\end{enumerate}
\vspace{1mm}
Then $\log(f(h))$ is constant and, consequently, $f(h)$ is constant. If, in addition, $f$ is not locally constant on $I$, then neither alternative can occur.}
\end{quote}

In the last section of this paper, it is shown how the previous uniqueness results provide new Calabi–Bernstein type results in the nonparametric case. First, in Corollary \ref{cor:corrected-graph-equation} we prove:
\begin{quote}{\it 
For a graph $\Sigma_u$ in a GRW spacetime $\Mbar$, where $u$ is a smooth function defined on an open subdomain $\Omega$ of $F$ and satisfying $u(\Omega)\subset I$, the conditions to be spacelike and 
$\langle A\nabla^{g_u}u,\nabla^{g_u}u\rangle_{g_u}=0$, where $\langle \, ,\,\rangle_{g_u}$ is the induced metric on $\Sigma_u$ \eqref{eq:induced-graph-metric}
and $\nabla^{g_u}u$ the corresponding gradient of $u$, 
are equivalent to

\vspace{1mm}

\begin{align}\tag{E}
\begin{cases}\displaystyle
  \Delta^F_\infty u
  +\frac{f'(u)}{f(u)}\big(f(u)^2-2|Du|_F^2\big)|Du|_F^2=0, \\[4mm]
  |Du|_F<f(u).
\end{cases}
\end{align} }
\end{quote}

In the case $\Omega=F$, a solution $u$ of (E), and hence $\Sigma_u$, is said entire. It should be noted that although $(F,\langle \, ,\, \rangle_F)$ is assumed complete, the induced metric \eqref{eq:induced-graph-metric} on $\Sigma_u$ is not complete, in general. Using a well-known sufficient condition for the completeness of $\Sigma_u$ where the constraint in (E) is changed to another stronger assumption, Proposition \ref{conditions-for-completeness}, we obtain in Theorem \ref{thm:graph-rigidity},

\begin{quote}{\it 
Let $f: I \to \R$ be a positive smooth function. Assume $f$ is not locally constant, has $\inf_I f>0$, and satisfies  $(\log f)'' \leq 0$. The only entire solutions $u\in C^\infty(F)$, $u(F)\subset I$, of the differential equation
\begin{align}\tag{E*}
\begin{cases}\displaystyle
            \Delta^F_\infty u
  +\frac{f'(u)}{f(u)}\big(f(u)^2-2|Du|_F^2\big)|Du|_F^2=0,\\[4mm]
            |Du|_F<\lambda f(u),
\end{cases}
\end{align}
where $\lambda$ is a constant with $0<\lambda<1$, on a complete Riemannian manifold $(F,\langle \, ,\, \rangle_F)$ are the constants $u=u_0\in I$.}
\end{quote} 

\vspace{1mm}

Among our other uniqueness results is Theorem \ref{thm:log-graph-rigidity},

\vspace{1mm}

\begin{quote}{\it 
Let $f:I\to\mathbb{R}$ be a positive smooth function. Assume $f$ is not locally constant and $\inf_I f>0$. The only entire solutions $u\in C^\infty(F)$, $u(F)\subset I$, of the differential equation \rm{(E*)}, that satisfy  $(\log f)''(u)\leq f'(u)^2/f(u)^2$
on a complete Riemannian manifold $(F,\langle \, ,\, \rangle_F)$ are the constants $u=u_0\in I$. }
\end{quote}

\vspace{1mm}
 
The remainder of the paper is organized as follows. In \Cref{sec:preliminaries}, we recall the basic notation and properties of GRW spacetimes, the height function, and the second fundamental form of a spacelike hypersurface. We also recall the Liouville-type property for the infinity-Laplacian obtained in \cite[p. 552, Eq. (1)]{AraujoMariPessoa2021} and compute the infinity-Laplacian of the height function in Proposition \ref{prop:infty-height}. In \Cref{sec:height}, we compute the infinity-Laplacian of three compositions of the height function: $f(h)$, a primitive $P(h)=\int^h f$, and $\log(f(h))$; see Propositions \ref{prop:infty-composed}, \ref{prop:primitive-formula}, and \ref{prop:logarithmic-height}, respectively. In \Cref{sec:main}, we prove the main global slice-rigidity theorems, two of which were stated above, with particular attention to the totally trapped case and the role of the TCC. Finally, in \Cref{sec:graphs}, we derive several Calabi--Bernstein-type uniqueness results for the differential equation (E*) from the rigidity theorems established in the preceding section, including the two results stated above.

\section{Preliminaries}\label{sec:preliminaries}

Let $(F,\langle\,,\,\rangle_F)$ be a connected\footnote{Throughout this paper, all manifolds are assumed to be connected.} $n(\geq 2)$-dimensional Riemannian manifold, let $I\subset\R$ be an open interval, and let $f:I\to(0,\infty)$ be a smooth function. We consider the generalized Robertson--Walker (GRW) spacetime
\begin{equation}\label{eq:grw-metric}
  \Mbar=I\times_f F,
  \qquad
  \langle\, ,\,\rangle =-\pi_I^*dt^2+f(\pi_I)^2\pi_F^*\langle\, ,\,\rangle_F,
\end{equation}
where $\pi_I:I\times F \to I$ and $\pi_F:I\times F \to F$ are the natural projections.
In other words, the $(n+1)$-dimensional spacetime $\Mbar$ is a warped product in the sense of \cite{ONeill1983}, with 
base $(I,-dt^2)$, fiber $(F,\langle\, ,\,\rangle_F)$ and warping function $f$. 
The coordinate vector field $\partial_t:=\partial/\partial_t$ is unit timelike, and the spacetime is time-oriented by $\partial_t$. 

\vspace{1mm}

On $\Mbar$ there exists a distinguished vector field, namely $f(\pi_I)\partial_t$, which is obviously timelike and satisfies 
\begin{equation}\label{eq:conformal}
\nbar_X\big(f(\pi_I)\partial_t\big)=f'(\pi_I)X,
\end{equation}
for any $X\in\mathfrak{X}(\Mbar)$, where $\nbar$ is the Levi-Civita connection of $\Mbar$. Thus, $f(\pi_I)\partial_t$ is conformal with $\mathcal{L}_{f(\pi_I)\partial_t}\langle\, ,\,\rangle = 2f'(\pi_I)\langle\, ,\,\rangle$
and its metrically equivalent $1$-form is closed. From (\ref{eq:conformal}) we get,
\begin{equation}\label{eq:divergence}
\overline{\text{div}}(\partial_t)=n\,\frac{f'(\pi_I)}{f(\pi_I)},
\end{equation}
where $\overline{\text{div}}$ denotes the divergence operator in $\Mbar$.

\vspace{1mm}

Let $\psi: M^n\longrightarrow \Mbar$ be a (connected) spacelike hypersurface in $\Mbar$, that is, $\psi$ is an immersion and the induced metric on $M$ via $\psi$ is Riemannian.  As usual, $\langle\,,\,\rangle$ also denotes the induced Riemannian metric on $M$.  Consider the unit normal vector field $N$ on $M$ chosen so that
\begin{equation}\label{eq:normal-orientation}
  \langle N,\partial_t\rangle<0,
\end{equation}
that is, future-pointing in the time orientation of $\Mbar$.

\vspace{1mm}

In what follows, we use the following sign convention for the second fundamental form, which is convenient in the literature on trapped submanifolds; see, for instance, \cite{AliasCanovasColares2017,AndradeCruzFigueiraLima2025}. The Weingarten formula for $\psi$ with respect to $N$ is written as
\begin{equation}\label{eq:weingarten}
 \nbar_X N=AX,
 \end{equation}
 for all $X\in\mathfrak{X}(M)$.
Moreover, the second fundamental form $\alpha$ is given by
\begin{equation}\label{eq:alpha-gr-convention}
  \alpha(X,Y)=-(\nbar_XY)^\perp,
\end{equation}
 for all $X,Y\in\mathfrak{X}(M)$. Thus, we have
\begin{equation}\label{eq:alpha-convention}
  \nbar_XY=\nabla_XY-\alpha(X,Y),
  \qquad
  \langle AX,Y\rangle=\langle\alpha(X,Y),N \rangle,
\end{equation}
where $\nabla$ is the Levi-Civita connection of $M$. Thus, we have
\begin{equation}\label{eq:alpha-hypersurface}
  \alpha(X,Y)=-\langle AX,Y\rangle N,
\end{equation}
and
\begin{equation}\label{eq:alpha-angle-relation}
  \langle \alpha(X,Y),\partial_t\rangle
  =-\langle N,\partial_t\rangle\langle AX,Y\rangle.
\end{equation}

Hence, for every nonzero tangent vector $X$, the normal vector $\alpha(X,X)$ is timelike and future-directed (resp. timelike and past-directed) precisely when $\langle AX, X\rangle<0$ (resp. $\langle AX, X\rangle>0$). Following the terminology of \cite{Senovilla2007,AndradeCruzFigueiraLima2025}, we say that $\psi$ is \emph{totally future trapped} (resp. \emph{totally past trapped}) if, for every point $p\in M$ and every nonzero vector $X\in T_pM$, the vector $\alpha(X,X)$ is timelike and future-directed (resp. timelike and past-directed). Since $\partial_t$ is future-directed, these conditions are equivalently written as
\begin{equation}\label{eq:strict-trapping-signs}
  \langle\alpha(X,X),\partial_t\rangle<0
  \qquad
  \left(\text{resp. }\langle\alpha(X,X),\partial_t\rangle>0\right)
\end{equation}
for every nonzero tangent vector $X$. In particular, the zero normal vector is not allowed by either definition.

\vspace{1mm}

When the strict trapping condition is evaluated at the height gradient, one obtains a weak inequality on all of $M$, because $\nabla h$ may vanish. More precisely, a totally future trapped (resp. totally past trapped) hypersurface satisfies
\begin{equation}\label{eq:trapping-gradient-signs}
  \langle\alpha(\nabla h,\nabla h),\partial_t\rangle\leq 0
  \qquad
  \left(\text{resp. }\langle\alpha(\nabla h,\nabla h),\partial_t\rangle\geq 0\right),
\end{equation}
and equality holds at a point if and only if $\nabla h$ vanishes there.

\begin{remark}\label{re:det(A)}
A spacelike hypersurface $\psi: M \rightarrow \Mbar$ is totally future trapped (resp. totally past trapped) if and only if the following symmetric tensor field (resp. the opposite of the following tensor field) $(X,Y) \mapsto -\langle AX, Y\rangle$, defined from the second fundamental form, is a Riemannian metric on $M^n$. In particular, we have det$(A)\neq 0$ everywhere on $M$. The converse is not true in general. In fact, if the Weingarten endomorphism satisfies this condition, then the metric previously defined on $M$ is non-degenerate. 
Because $M$ is connected, its index $s$ is constant. Thus, we have a semi-Riemannian metric on $M$. 
If $0<s<n=\text{dim}(M)$, then $M$ admits a distribution of dimension $s$ (and a complementary distribution of dimension $n-s$) \cite[Section 3.2]{GHV}. This structural property imposes topological restrictions on $M^n$ if it is assumed to be compact. In the special case $n=2m$ and $M$ is homeomorphic to $\mathbb{S}^{2m}$ we know that $M$ admits no non-trivial distribution \cite[Chap. 2]{Hatcher}. Therefore,  any spacelike hypersurface $M$ with the topology of $\mathbb{S}^{2m}$ satisfying $\det(A) \neq 0$ everywhere must be either totally future trapped or totally past trapped.
\end{remark}

The height function on the spacelike hypersurface $M$ in the GRW spacetime $\Mbar$ is defined by
\begin{equation}\label{eq:height}
  h=\pi_I\circ\psi.
\end{equation}
A spacelike hypersurface is said to be a spacelike slice when $h$ is constant \cite{AliasRomeroSanchez1995}, that is, if and only if $\psi(M)\subset\{t_0\}\times F$, for some $t_0\in I$.

\begin{remark}
For any spacelike hypersurface $\psi: M \to \Mbar$, in a GRW spacetime $\Mbar$, we know that $\pi_F\circ \psi: M \to F$ is a local diffeomorphism \cite[Sec. 3]{AliasRomeroSanchez1995}. Therefore, the projection on the fiber is always an open map and, in particular, $(\pi_F\circ \psi)(M)$ is an open subset of $F$. Moreover, 
in the particular case that $\psi$ is a spacelike slice, we can go further. Indeed, we have from (\ref{eq:grw-metric})
$$
\langle d(\pi_F\circ \psi)(v_1),d(\pi_F\circ \psi)(v_2)\rangle_F=f(t_0)^{-2}\langle v_1,v_2\rangle,
$$
for any $v_1,v_2\in T_xM$ and any $x\in M$, for each $t_0\in I$, i.e., $\pi_F\circ \psi$ is a local homothety. Consequently, if $M$ is also assumed to be complete, then $(\pi_F\circ \psi)(M)$ is also a closed subset of $F$. Indeed, a short geodesic segment joining a point of $(\pi_F\circ \psi)(M)$ to a point in its closure can be lifted locally to $M$; since the lifted segment has finite length, completeness provides an endpoint whose image is the prescribed limit point. By the connectedness of $F$, this gives $(\pi_F\circ \psi)(M)=F$.
\end{remark}

\vspace{1mm}

The tangential component $\partial_t^\top \in \mathfrak{X} (M)$ of $\partial_t$ along $\psi$ is
\begin{equation}\label{eq:decomposition}
  \partial_t^\top=\partial_t+\langle \partial_t,N\rangle N.
\end{equation}
For every $X\in \mathfrak{X}(M)$, we have
\begin{equation*}
  X(h)=-\langle \partial_t,X\rangle,
\end{equation*}
and therefore, the gradient of $h$ on $M$ is given by
\begin{equation}\label{eq:grad-height}
  \nabla h=-\partial_t^\top.
\end{equation}
Thus,
\begin{equation}\label{eq:height-angle}
  |\nabla h|^2=\langle \partial_t,N \rangle^2-1.
\end{equation}
Equivalently, if $\varphi$ denotes the hyperbolic angle between $\partial_t$ and $N$, then $\langle \partial_t,N \rangle=-\cosh\varphi$ and we can rewrite (\ref{eq:height-angle}) as follows 
\begin{equation}\label{eq:height-angle2}
|\nabla h|^2=\sinh^2\varphi.
\end{equation}
Thus, a spacelike hypersurface is a spacelike slice if and only if the hyperbolic angle $\varphi$ vanishes identically.

\vspace{1mm}

Let $(P, g)$ be an $n$-dimensional Riemannian manifold. 
The infinity Laplacian operator on $(P, g)$ is defined for each $v \in C^2(P)$ by
\begin{equation}\label{eq:infty-def}
\Delta_{\infty}v=\Hess v\,(\nabla v,\nabla v)
  =g(\nabla_{\nabla v}\nabla v,\nabla v),
\end{equation}
where $\nabla$ denotes the Levi-Civita connection of $g$ and $\text{Hess}\, v$ 
is the Riemannian Hessian of $v$.  The infinity Laplacian naturally arises in geometric analysis as the Euler--Lagrange operator for absolutely minimizing Lipschitz extensions on metric-measure spaces and Riemannian manifolds (see, e.g., \cite{Jensen1993,Juutinen2001}). 
Alternatively, $\Delta_\infty$ can be understood as the asymptotic limit 
as $p \to \infty$ of the normalized $p$-Laplacian operators 
$$
\Delta_p^N u := \frac{1}{p-2}\, |\nabla u|^{2-p} \,\text{div}_g\big(|\nabla u|^{p-2}\nabla u\big).
$$
Although $\Delta_{\infty}$ is degenerate and nonlinear, it shows an important Liouville-type property, which follows as a consequence of the potential-theoretic characterization of completeness through the infinity-Laplacian due to Araujo, Mari and Pessoa \cite[p. 552, Eq. (1)]{AraujoMariPessoa2021}. 

\begin{lemma}\label{lem:liouville}
Let $(P,g)$ be a complete Riemannian manifold. If $v\in C^2(P)$ satisfies $\Delta_\infty v\geq 0$ (resp. $\Delta_\infty v\leq 0$) and is bounded from above (resp. from below), then $v$ must be constant. 
\end{lemma}

The next result shows the basic formula used throughout the paper. It can be seen as the spacelike hypersurface counterpart of the corresponding infinity-Laplacian computation used for codimension two spacelike submanifolds in \cite{AndradeCruzFigueiraLima2025}.

\begin{proposition}\label{prop:infty-height}
Let $\psi: M^n\to \Mbar$ be a spacelike hypersurface in a GRW spacetime $\Mbar$. Then,  the height function $h=\pi_I\circ\psi$ satisfies 
\begin{equation}\label{eq:infty-height}
  \Delta_\infty h
  =-\frac{f'(h)}{f(h)}\big(1+|\nabla h|^2\big)|\nabla h|^2
  -\langle N,\partial_t\rangle\langle A\nabla h,\nabla h\rangle.
\end{equation}
\end{proposition}

\begin{proof}
From (\ref{eq:conformal}), we have
\begin{equation}\label{eq:connection-partialt}
  \nbar_X\partial_t=\frac{f'(h)}{f(h)}\left(X+\langle X,\partial_t\rangle\partial_t\right)
\end{equation}
for every $X\in \mathfrak{X}(M)$. From \eqref{eq:decomposition}, using \eqref{eq:weingarten} and \eqref{eq:alpha-convention}, formula \eqref{eq:connection-partialt} gives 
$$
\nabla_X\partial^{\top}_t=\frac{f'(h)}{f(h)}\left(X+\langle X,\partial_t\rangle\partial^{\top}_t\right) + \langle \partial_t,N\rangle AX.
$$
Taking into account (\ref{eq:grad-height}), previous formula can be rewritten
\begin{equation}\label{eq:hess-height-1}
\nabla_{X}\nabla h=-\frac{f'(h)}{f(h)}\left(X + \langle X,\partial_t\rangle\partial^{\top}_t\right) - \langle \partial_t,N\rangle AX,
\end{equation}
and, therefore,
\begin{equation}\label{eq:hess-height}
  \Hess h(X,Y)
  =-\frac{f'(h)}{f(h)}
  \left(\langle X,Y\rangle+\langle X,\partial_t\rangle\langle Y,\partial_t\rangle\right)-\langle N,\partial_t\rangle\langle AX,Y\rangle.
\end{equation}
Finally, taking $X=Y=\nabla h$ and using again \eqref{eq:grad-height} in the previous formula, we obtain \eqref{eq:infty-height}, ending the proof.
\end{proof}

\begin{remark}\label{rem:height-alpha-form}
Using \eqref{eq:alpha-angle-relation}, identity \eqref{eq:infty-height} can be rewritten in terms of the second fundamental form as
\begin{equation}\label{eq:infty-height-alpha}
  \Delta_\infty h
  = -\frac{f'(h)}{f(h)}\big(1+|\nabla h|^2\big)|\nabla h|^2+\langle \alpha(\nabla h,\nabla h),\partial_t\rangle.
\end{equation}
\end{remark}

We end this section with two simple technical results to be used later.

\vspace{1mm}

\begin{lemma}\label{lem:chain-infty-b}
For every smooth function $\rho : I\to\R$ and every smooth $u: P\to I$, the infinity-Laplacian of the composition  $\rho (u)=\rho \circ u$ satisfies
\begin{equation}\label{eq:chain-infty}
  \Delta_\infty(\rho(u))
  =\rho'(u)^2\left\{\rho'(u)\Delta_\infty u+\rho''(u)|\nabla u|^4\right\}.
\end{equation}
\end{lemma}

\begin{proof}
This easily follows from
\begin{equation*}
  \nabla(\rho(u))=\rho'(u)\nabla u,
  \qquad
  \Hess(\rho(u))=\rho''(u)du\otimes du+\rho'(u)\Hess u.
\end{equation*}

\vspace{-6mm}

\end{proof}

A smooth function defined on an open interval $I\subset\R$, $f: I\to\R$, is said to be not locally constant if there is no nonempty open sub-interval $J\subset I$ such that $f|_J$ is constant. 

\begin{lemma}\label{lem:composition-constant}
If a smooth function $f: I\to\R$ is not locally constant and $f\circ u$ is constant for a smooth function $u: P\to I$ on a connected manifold $P$, then $u$ must be constant.
\end{lemma}

\begin{proof}
Since $P$ is connected, $u(P)$ is a connected subset of $\R$. If $u(P)$ contains more than one point, then it contains a nonempty open subinterval of $I$. Since $f$ is constant on $u(P)$, the restriction of $f$ to such an interval would be constant, which contradicts the assumption. Hence, $u(P)$ is a single point.
\end{proof}

\begin{remark}\label{rem:geometrical-meaning}
If the warping function $f: I \to \mathbb{R}$ of a GRW spacetime $\Mbar$ is not locally constant, then there exists no open subinterval $J$ of $I$ such that $J\times_{f|_{J}}F$ is a semi-Riemannian product. 
\end{remark}
\section{Compositions of the height function}\label{sec:height}

In this section, we develop the three basic height mechanisms used throughout the paper. They arise from the infinity-Laplacian of the warping composition $f(h)$, of a primitive composition $P(h)=\int^h f$, and of the logarithmic composition $\log(f(h))$. These identities will later be combined with completeness and sign conditions to obtain slice rigidity.

\subsection{The warping composition}\label{subsec:warping-composition}

We first record the identity associated with the direct composition of the height function with the warping function.

\begin{proposition}\label{prop:infty-composed}
Let $\psi: M^n\to \Mbar$ be a spacelike hypersurface in a GRW spacetime $\Mbar$. Then, the function $f(h)=f\circ \pi_I\circ\psi$ satisfies 
\begin{equation}\label{eq:infty-fheight}
    \begin{aligned}
         \Delta_\infty(f(h))
  ={}&f'(h)^2\Big\{
  -\frac{f'(h)^2}{f(h)}|\nabla h|^2
  +f(h)(\log f)''(h)|\nabla h|^4  \\
  &\hspace{-0.5cm}
  +f'(h)\langle \alpha(\nabla h,\nabla h),\partial_t\rangle
  \Big\}.
    \end{aligned}
\end{equation}
\end{proposition}
\begin{proof}
Taking $\rho=f$ and $u=h$ in Lemma \ref{lem:chain-infty-b}, and then using Proposition \ref{prop:infty-height}, we obtain
\begin{align*}
  \Delta_\infty(f(h))
  ={}&f'(h)^2\Big\{
  -\frac{f'(h)^2}{f(h)}\big(1+|\nabla h|^2\big)|\nabla h|^2
  +f''(h)|\nabla h|^4 \\
  &\hspace{-0.5cm}
  +f'(h)\langle \alpha(\nabla h,\nabla h),\partial_t\rangle
  \Big\}.
\end{align*}
Since
\begin{equation*}
  f''-\frac{(f')^2}{f}=f(\log f)'',
\end{equation*}
the two terms of order $|\nabla h|^4$ combine and give \eqref{eq:infty-fheight}.
\end{proof}

\subsection{The primitive test function}\label{subsec:primitive}

There is another composition that is useful when the warping function is monotone. Fix $t_*\in I$ and define
\begin{equation}\label{eq:primitive-def}
  P(t)=\int_{t_*}^{t}f(s)\,ds,
\end{equation}
that is, $P$ is a primitive of $f$. Then
\begin{equation*}
  P'(t)=f(t)>0,
  \qquad
  P''(t)=f'(t),
\end{equation*}
in particular, $P$ is strictly increasing on $I$.

\begin{proposition}\label{prop:primitive-formula}
Let $\psi: M^n\to \Mbar$ be a spacelike hypersurface in a GRW spacetime $\Mbar$, and let $P$ be a primitive of $f$. Then
\begin{equation}\label{eq:primitive-formula}
  \Delta_\infty(P(h))
  =f(h)^2\left[f(h)\langle \alpha(\nabla h,\nabla h),\partial_t\rangle
  -f'(h)|\nabla h|^2\right].
\end{equation}
\end{proposition}

\begin{proof}
Using the chain rule \eqref{eq:chain-infty} with $\rho=P$ and $u=h$, we get
\begin{equation*}
  \Delta_\infty(P(h))
  =f(h)^2\left\{f(h)\Delta_\infty h+f'(h)|\nabla h|^4\right\}.
\end{equation*}
Now, from Proposition \ref{prop:infty-height},
\begin{align*}
  \Delta_\infty(P(h))
  ={}&f(h)^2\Big\{
  f(h)\langle \alpha(\nabla h,\nabla h),\partial_t\rangle
  -f'(h)(1+|\nabla h|^2)|\nabla h|^2 
  +f'(h)|\nabla h|^4
  \Big\}.
\end{align*}
The two quartic terms cancel, and this gives \eqref{eq:primitive-formula}.
\end{proof}

\subsection{The logarithmic warping composition}\label{subsec:logarithmic-composition}

We now consider the logarithmic counterpart of Proposition \ref{prop:infty-composed}. This computation is useful because it can be written entirely in terms of the vector-valued second fundamental form $\alpha$. Set
\begin{equation}\label{eq:eta-def}
  \eta(t)=(\log f)'(t)=\frac{f'(t)}{f(t)}.
\end{equation}
Then
\begin{equation}\label{eq:eta-prime-def}
  \eta'(t)=(\log f)''(t)
  =\frac{f''(t)}{f(t)}-\frac{f'(t)^2}{f(t)^2}.
\end{equation}

\begin{proposition}\label{prop:logarithmic-height}
Let $\psi: M^n\to \Mbar$ be a spacelike hypersurface in a GRW spacetime $\Mbar$, and let $\log(f(h))=\log \circ f \circ \pi_I\circ\psi$. Then
\begin{equation}\label{eq:log-alpha-formula}
\begin{aligned}
  \Delta_\infty\log(f(h))
  ={}&
  \eta(h)^2\Big[
  \eta(h)\langle \alpha(\nabla h,\nabla h),\partial_t\rangle
  -\eta(h)^2|\nabla h|^2 \\
  &\hspace{-0.5cm}
  +\big(\eta'(h)-\eta(h)^2\big)|\nabla h|^4
  \Big].
\end{aligned}
\end{equation}
\end{proposition}

\begin{proof}
Let $\rho(t)=\log f(t)$. Then $\rho'=\eta$ and $\rho''=\eta'$. Applying the chain rule \eqref{eq:chain-infty} with $u=h$ gives
\begin{equation}\label{eq:log-chain-height}
  \Delta_\infty(\rho(h))
  =\eta(h)^2\left\{\eta(h)\Delta_\infty h+\eta'(h)|\nabla h|^4\right\}.
\end{equation}
Substituting the height identity \eqref{eq:infty-height-alpha} we get \eqref{eq:log-alpha-formula}. 
\end{proof}
\begin{remark}
 Equivalently, in terms of $f$, we have
\begin{equation}\label{eq:log-alpha-formula-f}
\begin{aligned}
  \Delta_\infty\log(f(h))
  ={}&
  \left(\frac{f'(h)}{f(h)}\right)^2\,
  \Big[
  \left(\frac{f'(h)}{f(h)}\right)
  \langle \alpha(\nabla h,\nabla h),\partial_t\rangle
  \\[1mm]
  &\hspace{-0.5cm}
  -\left(\frac{f'(h)}{f(h)}\right)^2\,|\nabla h|^2
  +
  \left(
  \frac{f''(h)}{f(h)}-2\frac{f'(h)^2}{f(h)^2}
  \right)\,|\nabla h|^4
  \Big],
\end{aligned}
\end{equation}
that follows from 
\begin{equation*}
  \eta'(h)-\eta(h)^2
  =\frac{f''(h)}{f(h)}-2\frac{f'(h)^2}{f(h)^2}.
\end{equation*}   
\end{remark}

\begin{remark}\label{rem:log-vs-primitive}
The primitive identity and the logarithmic identity should be read as complementary mechanisms. For the primitive $P(h)$, the quartic gradient terms cancel exactly in \eqref{eq:primitive-formula}. For the logarithmic composition, no such cancellation occurs; instead, the sign of the quartic term is controlled by the inequality 
\begin{equation}\label{eq:log-warping-condition}
  \eta'(h)\leq \eta(h)^2,
\end{equation}
which is equivalent to
\begin{equation}\label{eq:log-warping-condition2}
(\log f)''(h)\leq \frac{f'(h)^2}{f(h)^2}\,.
\end{equation}

We end this section by exposing two energy conditions on spacetime, especially on GRW spacetime, and their physical interpretation. Recall now that a spacetime satisfies the timelike convergence condition (TCC) (resp. the null convergence condition (NCC)) if
\begin{equation*}
  \Ricbar(Z,Z)\geq 0,
\end{equation*}
for every timelike tangent vector $Z$ (resp. for every null tangent vector $Z$) where $\Ricbar$ is the Ricci tensor of the spacetime. 

\vspace{1mm}

The TCC has a natural physical interpretation. In fact, the TCC on a spacetime may be interpreted as expressing that gravity is, on average, attractive. In four-dimensional General Relativity, under the Einstein field equation with vanishing cosmological constant, the TCC is related to suitable energy conditions imposed on the stress-energy tensor; see \cite{SachsWu1977}.  
A continuity argument at each point in spacetime shows that the TCC implies the NCC, and the converse is not true in general.

\vspace{1mm }

A GRW spacetime  $\Mbar=I\times_f F$, satisfies the NCC if and only if
\begin{equation}\label{eq:ncc-grw}
\Ric^F(Y,Y)\geq (n-1)f(t)^2(\log f)''(t)\langle Y,Y\rangle_F
\end{equation}
for every $Y\in T_xF$, $(t,x)\in \Mbar$, where $\Ric^F$ is the Ricci tensor of $F$. The GRW spacetime $\Mbar$ satisfies the TCC if and only if it satisfies the NCC and $f''\leq 0$ (see, for instance, \cite{AledoRubioSalamanca2017,AliasRomeroSanchez1995}).

\vspace{1mm}

It should be pointed out that the NCC alone does not imply $(\log f)''\leq 0$, i.e., the convexity of $-\log f$ in general. The implication holds if we assume $\Ric^F$ is non-positive definite. Of course, the TCC does imply that $-\log f$ is convex.

\vspace{1mm}

On the other hand, from (\ref{eq:divergence}) we have that $f'<0$ (resp. $f'>0$) may be interpreted as saying that the observers in the reference frame $\partial_t$ are on average coming together (resp. spreading apart).  If we assume $-\log f$ convex then 
$$
\frac{d}{ds}\big(\overline{\text{div}}(\partial_t)\circ \gamma \big)(s)\leq 0,
$$
for any observer $\gamma$ in $\partial_t$. If in addition we assume there is a proper time $s_0$ such that $\overline{\text{div}}(\partial_t)_{\gamma (s_0)}<0$, then $\overline{\text{div}}(\partial_t)_{\gamma (s)}<0$ for any $s>s_0$. 
Thus, the convexity of $-\log f$ does not force expansion, but constrains the expansion rate to be non-increasing along the observers of the reference frame $\partial_t$.
\end{remark}

\section{Rigidity criteria for spacelike slices}\label{sec:main}

We now apply the results of the previous section to obtain slice rigidity theorems for certain complete spacelike hypersurfaces. We divide the discussion according to the three test functions used: the warping composition $f(h)$, the primitive composition $P(h)$, and the logarithmic composition $\log(f(h))$.

\subsection{Rigidity results for the warping criterion}\label{subsec:warping-criterion}

The first criterion is obtained from (\ref{eq:infty-fheight}).

\begin{theorem}\label{thm:main-rigidity}
Let $\psi: M^n\to \Mbar$ be a complete spacelike hypersurface in a GRW spacetime $\Mbar$, and let $h=\pi_I\circ\psi$. Assume 
\begin{equation}\label{eq:convex-assumption}
  (\log f)''(h)\leq 0
  \qquad\text{on} \quad M,
\end{equation}
and
\begin{equation}\label{eq:directional-sign}
  f'(h)\langle \alpha(\nabla h,\nabla h),\partial_t\rangle\leq0
  \qquad\text{on} \quad M.
\end{equation}
Then $f(h)$ is constant. If, in addition, $f$ is not locally constant on $I$, then $h$ is a constant $t_0\in I$. Consequently, we have $\psi(M)=\{t_0\}\times F$.
\end{theorem}

\begin{proof}
Since $f>0$, the function $f(h)$ is bounded from below by $0$. By Proposition \ref{prop:infty-composed}, assumptions \eqref{eq:convex-assumption} and \eqref{eq:directional-sign} imply
\begin{equation*}
  \Delta_\infty(f(h))\leq0.
\end{equation*}
Because $M$ is complete, Lemma \ref{lem:liouville} gives that $f(h)$ is constant. If $f$ is not locally constant, Lemma \ref{lem:composition-constant} implies that $h=t_0$ is constant for some $t_0\in I$. This ends the proof.
\end{proof}

The assumption (\ref{eq:directional-sign}) clearly is satisfied if the vanishing of the directional second fundamental term holds.

\begin{corollary}\label{cor:directional-zero}
Let $\psi: M^n\to \Mbar$ be a complete spacelike hypersurface in a GRW spacetime $\Mbar$, and let $h=\pi_I\circ\psi$. Suppose that
\begin{equation*}
  (\log f)''(h)\leq0
  \qquad\text{and}\qquad
  \langle A\nabla h,\nabla h\rangle=0
\end{equation*}
on $M$. Then $f(h)$ is constant. Moreover, if $f$ is not locally constant, then $h$ is constant $t_0\in I$ and we have $\psi(M)=\{t_0\}\times F$.
\end{corollary}

\begin{remark}\label{rem:directional-condition}
The hypothesis $\langle A\nabla h,\nabla h\rangle=0$ is directional. It should not be replaced by $\langle AX, X\rangle=0$ for every tangent vector $X$, which means $A=0$. Indeed, a spacelike slice $\{t_0\}\times F$ has
\begin{equation*}
  A=\frac{f'(t_0)}{f(t_0)}I,
\end{equation*}
with the convention \eqref{eq:weingarten}. Hence $A$ is not identically zero unless $f'(t_0)=0$. However, on a slice, one has $\nabla h=0$, and therefore the directional condition is automatically satisfied.
\end{remark}
\begin{counter-example}
Consider the $(n+1)$-dimensional De Sitter spacetime of sectional curvature $1$, $\mathbb{S}^{n+1}_1=\{(x_0,x_1,...,x_{n+1})\in\mathbb{L}^{n+2}\, : \, -x_0^2+x_1^2+...+x_{n+1}^2=1\}$. It is well-known that its totally geodesic compact spacelike hypersurfaces are $\mathbb{S}_1^{n+1}\cap H_v$, where $v$ is a timelike vector of $\mathbb{L}^{n+2}$ and $H_v$ is the linear hyperplane of $\mathbb{L}^{n+2}$ with normal vector $v$. In the representation of $\mathbb{S}^{n+1}_1$ as a GRW spacetime with base $(\mathbb{R},-dt^2)$, fiber $(\mathbb{S}^n,g_{\mathbb{S}^n})$ and warping function $f(t)=\cosh(t)$, the only totally geodesic slice is $t=0$. Any other totally geodesic compact spacelike hypersurface $M$ clearly satisfies $\langle A\nabla h,\nabla h\rangle=0$; however, $h$ is not constant. Note that inequality (\ref{eq:directional-sign}) holds trivially; however, \eqref{eq:convex-assumption} does not hold. Thus, this assumption cannot be dropped in previous results. 
\end{counter-example}
The sign assumption (\ref{eq:directional-sign}) in Theorem \ref{thm:main-rigidity} is automatically satisfied when the causal orientation of $\alpha$ is compatible with the monotonicity of the warping function.

To illustrate the applications of Corollary \ref{cor:directional-zero}, consider the $(n+1)$-dimensional Steady State spacetime that can be viewed as the GRW spacetime with base $(\mathbb{R},-dt^2)$, fiber the Euclidean space $(\mathbb{R}^n,g^o)$ and warping function $f(t)=e^t$ (see, for instance, \cite{AlbujerAlias2009}). Note that in this case we have $(\log f)''=0$. Thus, Corollary \ref{cor:directional-zero} specializes to the following uniqueness result  
\begin{corollary}\label{Steady_State}
The only complete spacelike hypersurfaces in the Steady State spacetime satisfying 
\begin{equation*}
\langle A\nabla h,\nabla h\rangle=0,
\end{equation*}
are the spacelike slices $\{t_0\}\times \mathbb{R}^n$.
\end{corollary}
\begin{corollary}\label{cor:totally-trapped-version}
Let $\Mbar=I\times_fF$ be a GRW spacetime whose warping function is not locally constant. There exists no complete totally future trapped spacelike hypersurface $\psi:M^n\to\Mbar$ satisfying
\begin{equation*}
  (\log f)''(h)\leq0
  \qquad\text{and}\qquad
  f'(h)\geq0
\end{equation*}
on $M$, where $h=\pi_I\circ\psi$. Likewise, there exists no complete totally past trapped spacelike hypersurface satisfying
\begin{equation*}
  (\log f)''(h)\leq0
  \qquad\text{and}\qquad
  f'(h)\leq0
\end{equation*}
on $M$.
\end{corollary}

\begin{proof}
Suppose first that such a totally future trapped hypersurface exists. By \eqref{eq:trapping-gradient-signs},
\begin{equation*}
  \langle\alpha(\nabla h,\nabla h),\partial_t\rangle\leq0,
\end{equation*}
and hence
\begin{equation*}
  f'(h)\langle\alpha(\nabla h,\nabla h),\partial_t\rangle\leq0.
\end{equation*}
Theorem \ref{thm:main-rigidity} implies that $h=t_0$ is constant and that $\psi(M)=\{t_0\}\times F$. On this slice, the future-pointing unit normal is $N=\partial_t$ and
\begin{equation*}
  A=\frac{f'(t_0)}{f(t_0)}I.
\end{equation*}
Consequently, for every nonzero tangent vector $X$, \eqref{eq:alpha-angle-relation} gives
\begin{equation*}
  \langle\alpha(X,X),\partial_t\rangle
  =\frac{f'(t_0)}{f(t_0)}|X|^2\geq0,
\end{equation*}
which contradicts the strict inequality in \eqref{eq:strict-trapping-signs} for a totally future trapped hypersurface.

\vspace{1mm}

The totally past trapped case is analogous. Indeed, Theorem \ref{thm:main-rigidity} again forces the hypersurface to be a slice, where
\begin{equation*}
  \langle\alpha(X,X),\partial_t\rangle
  =\frac{f'(t_0)}{f(t_0)}|X|^2\leq0,
\end{equation*}
contradicting the strict positive sign required by \eqref{eq:strict-trapping-signs}.
\end{proof}

The preceding results give the following application to the $(n+1)$-dimensional Einstein--De Sitter spacetime. Note that, for $n=3$, this is the usual four-dimensional Einstein--De Sitter spacetime (see, for instance, \cite{SachsWu1977}).

\begin{corollary}\label{cor:einstein-de-sitter}
Consider the GRW spacetime
\begin{equation*}
  \Mbar_{\mathrm{EdS}}
  =(0,\infty)\times_{t^{2/3}}\mathbb{R}^n,
  \qquad
  \langle\, ,\,\rangle=-dt^2+t^{4/3}g^o.
\end{equation*}
Then the following statements hold:
\begin{enumerate}
  \item[\rm(1)] The only complete spacelike hypersurfaces satisfying
  \begin{equation*}
    \langle A\nabla h,\nabla h\rangle=0
  \end{equation*}
  are the spacelike slices $\{t_0\}\times\mathbb{R}^n$, with $t_0>0$.

\vspace*{1mm}

  \item[\rm(2)] There exists no complete totally future trapped spacelike hypersurface in $\Mbar_{\mathrm{EdS}}$.
\end{enumerate}
\end{corollary}

\begin{proof}
For $f(t)=t^{2/3}$, one has
\begin{equation*}
  f'(t)=\frac{2}{3}t^{-1/3}>0,
  \qquad
  (\log f)''(t)=-\frac{2}{3t^2}<0.
\end{equation*}
Moreover, $f$ is strictly increasing. Assertion (1) follows from Corollary \ref{cor:directional-zero}, whereas assertion (2) follows from Corollary \ref{cor:totally-trapped-version}.
\end{proof}

\subsection{Rigidity results for the primitive composition criterion}\label{subsec:primitive-criterion}

The next result uses the primitive composition $P(h)$ and the cancellation displayed in Proposition \ref{prop:primitive-formula}.

\begin{theorem}\label{thm:primitive-monotonicity}
Let $\psi: M^n\to \Mbar$ be a complete spacelike hypersurface in a GRW spacetime $\Mbar$, and let $h=\pi_I\circ\psi$. Assume 
\begin{equation}\label{eq:primitive-directional-zero}
  \langle A\nabla h,\nabla h\rangle=0
\end{equation}
on $M$. If one of the following assertions holds:
\begin{enumerate}
  \item[\rm(1)] $f'(h)\geq0$ on $M$ and $P(h)$ is bounded from below, 
  \vspace*{1mm}
  \item[\rm(2)] $f'(h)\leq0$ on $M$ and $P(h)$ is bounded from above, 
\end{enumerate}
then, $h=t_0$ for $t_0\in I$  and $\psi(M)=\{t_0\}\times F$. 
\end{theorem}

\begin{proof}
By Proposition \ref{prop:primitive-formula} and formula \eqref{eq:primitive-directional-zero},
\begin{equation}\label{eq:primitive-zero-reduction}
  \Delta_\infty(P(h))=-f(h)^2f'(h)|\nabla h|^2.
\end{equation}

If $f'(h)\geq0$, then $\Delta_\infty(P(h))\leq0$. Since $P(h)$ is bounded from below and $M$ is complete, Lemma \ref{lem:liouville} implies that $P(h)$ is constant. Because $P$ is strictly increasing, $h$ is constant.

\vspace{1mm}

If $f'(h)\leq0$, then \eqref{eq:primitive-zero-reduction} gives $\Delta_\infty(P(h))\geq0$. Since $P(h)$ is bounded from above, Lemma \ref{lem:liouville} again implies that $P(h)$ is constant. Hence $h$ is constant.
\end{proof}

\begin{remark}\label{rem:primitive-boundedness}
The boundedness assumptions in Theorem \ref{thm:primitive-monotonicity} are naturally formulated in terms of $P(h)$. Since $P$ is increasing, $P(h)$ is bounded from below, for instance, whenever $h$ is bounded from below by some value $t_-\in I$. Similarly, $P(h)$ is bounded from above whenever $h$ is bounded from above by some value $t_+\in I$.
\end{remark}

\subsection{Rigidity results for the logarithmic criterion}\label{subsec:logarithmic-criterion}

For the logarithmic composition, the sign estimate depends on the milder condition
(\ref{eq:log-warping-condition}).
Thus, in this case, one does not need to assume the convexity of $-\log f$ (although if it holds, then inequality (\ref{eq:log-warping-condition2}) is satisfied). The timelike convergence condition enters the discussion simply as a geometric way of ensuring that this inequality holds.

\begin{theorem}\label{thm:log-trapped-rigidity}
Let $\psi: M^n\to\Mbar$ be a complete spacelike hypersurface in a GRW spacetime $\Mbar=I\times_f F$, and let $h=\pi_I\circ\psi$ be its height function. Set
\begin{equation*}
  \eta(t)=\frac{f'(t)}{f(t)}
\end{equation*}
and assume 
\begin{equation}\label{eq:log-warping-main}
  \eta'(h)\leq\eta(h)^2
  \qquad\text{on }M,
\end{equation}
and
\begin{equation}\label{eq:log-lower-bound-main}
  f(h)\geq c>0
  \qquad\text{on}\; M,
\end{equation}
for some constant $c>0$. Suppose that one of the following alternatives holds:
\begin{enumerate}
  \item[\rm(1)] $\psi$ is totally future trapped and $f'(h)\geq0$ on $M$;
\vspace*{1mm}
  \item[\rm(2)] $\psi$ is totally past trapped and $f'(h)\leq0$ on $M$.
\end{enumerate}
Then $\log(f(h))$ is constant on $M$, and consequently $f(h)$ is constant. If, in addition, $f$ is not locally constant on $I$, then neither alternative can occur. In other words, under this additional assumption, there exists no complete totally future trapped spacelike hypersurface satisfying alternative {\rm(1)}, and no complete totally past trapped spacelike hypersurface satisfying alternative {\rm(2)}.
\end{theorem}

\begin{proof}
Suppose first that $\psi$ is totally future trapped. By \eqref{eq:trapping-gradient-signs},
\begin{equation*}
  \langle\alpha(\nabla h,\nabla h),\partial_t\rangle\leq0.
\end{equation*}
Since $f'(h)\geq0$ and $f(h)>0$, we have $\eta(h)\geq0$, and hence
\begin{equation*}
  \eta(h)\langle\alpha(\nabla h,\nabla h),\partial_t\rangle\leq0.
\end{equation*}
If $\psi$ is totally past trapped, then
\begin{equation*}
  \langle\alpha(\nabla h,\nabla h),\partial_t\rangle\geq0.
\end{equation*}
In this case, $f'(h)\leq0$ implies $\eta(h)\leq0$, and therefore the same inequality holds:
\begin{equation*}
  \eta(h)\langle\alpha(\nabla h,\nabla h),\partial_t\rangle\leq0.
\end{equation*}
In either case, the preceding sign inequality and \eqref{eq:log-warping-main} imply that
\begin{equation*}
  \Delta_\infty\log(f(h))\leq0.
\end{equation*}
Moreover, \eqref{eq:log-lower-bound-main} yields $\log(f(h))\geq\log c$. Since $M$ is complete, Lemma \ref{lem:liouville} implies that $\log(f(h))$ is constant. Consequently, $f(h)$ is constant.

\vspace{1mm}

Assume now that $f$ is not locally constant on $I$. By Lemma \ref{lem:composition-constant}, the constancy of $f(h)$ implies that $h=t_0$ for some $t_0\in I$. Thus,
\begin{equation*}
  \psi(M)=\{t_0\}\times F.
\end{equation*}
On this slice, the future-pointing unit normal is $N=\partial_t$ and
\begin{equation*}
  A=\frac{f'(t_0)}{f(t_0)}I.
\end{equation*}
Hence, for every nonzero tangent vector $X$,
\begin{equation}\label{eq:log-slice-sign}
  \langle\alpha(X,X),\partial_t\rangle
  =\frac{f'(t_0)}{f(t_0)}|X|^2.
\end{equation}
Under alternative {\rm(1)}, the right-hand side of \eqref{eq:log-slice-sign} is nonnegative, contradicting the strict negative sign required for a totally future trapped hypersurface. On the other hand, under alternative {\rm(2)}, it is non-positive, contradicting the strict positive sign required for a totally past trapped hypersurface. Therefore, if $f$ is not locally constant, neither alternative can occur.
\end{proof}

\begin{corollary}\label{cor:tcc-log-rigidity}
Let $\psi: M^n\to \Mbar$ be a complete spacelike hypersurface in a GRW spacetime $\Mbar$. Assume $\Mbar$ satisfies the TCC and 
\begin{equation*}
  f(h)\geq c>0,
\end{equation*}
for some constant $c>0$. Suppose that one of the following alternatives holds:
\begin{enumerate}
  \item[\rm(1)] $\psi$ is totally future trapped and $f'(h)\geq0$ on $M$,

\vspace*{1mm}
  
  \item[\rm(2)] $\psi$ is totally past trapped and $f'(h)\leq0$ on $M$.
\end{enumerate}
Then $\log(f(h))$ is constant and, consequently, $f(h)$ is constant. If, in addition, $f$ is not locally constant on $I$, then neither alternative can occur.
\end{corollary}

\begin{proof}
Under the TCC we have $\Ricbar(\partial_t,\partial_t)\geq 0$, and therefore $f''\leq0$. Hence
\begin{equation*}
  \eta'-\eta^2
  =\frac{f''}{f}-2\frac{(f')^2}{f^2}
  \leq0,
\end{equation*}
thus (\ref{eq:log-warping-condition}) holds. The result follows from Theorem \ref{thm:log-trapped-rigidity}.
\end{proof}

\begin{remark}\label{rem:log-result-role}
 The estimate for $f(h)$ uses the sign of $(\log f)''$, whereas the logarithmic estimate uses the condition $\eta'\leq\eta^2$, or equivalently (\ref{eq:log-warping-condition2}). This allows the logarithmic criterion to cover situations not included in the hypothesis $(\log f)''\leq 0$. Under the TCC, the logarithmic warping inequality is automatic. The price paid is the lower bound $f(h)\geq c>0$, needed to apply the Liouville principle Lemma \ref{lem:liouville} to $\log(f(h))$.
\end{remark}

\begin{corollary}\label{cor:ncc-applications}
Let $\psi: M^n\to \Mbar$ be a complete spacelike hypersurface in a GRW spacetime $\Mbar$, satisfying
\begin{equation*}
  \langle A\nabla h,\nabla h\rangle=0.
\end{equation*}
Suppose that the warping function $f$ is not locally constant and one of the following alternatives holds: 
\begin{enumerate}
  \item[\rm(1)] $\Mbar$ satisfies the TCC,
\vspace*{1mm}
  \item[\rm(2)] $\Mbar$ satisfies the NCC and the Ricci tensor of the fiber $F$ is non-positive definite.
\end{enumerate}
Then $h$ is a constant $t_0\in I$ and $\psi(M)=\{t_0\}\times F$. 
\end{corollary}

\begin{proof}
If $\Mbar$ satisfies the TCC, then in particular the warping function satisfies $f''\leq0$, hence $(\log f)'' \leq 0$, which gives the conclusion under (1). If the NCC holds and $\Ric^F\leq0$, then \eqref{eq:ncc-grw} forces $(\log f)''\leq0$, proving the conclusion under (2). 
\end{proof}

\begin{remark}\label{rem:einstein-grw}
If $\Mbar$ is assumed to be Einstein, then it clearly satisfies the NCC. In the classification of Einstein GRW spacetimes, the cases with non-positive Ricci curvature on the fiber give natural examples where the previous corollary applies. This includes the corresponding cases in the table of Einstein GRW spacetimes given in \cite{AliasRomeroSanchez1997}.
\end{remark}

\section{Spacelike graphs and the non-parametric equation}\label{sec:graphs}

The graphical case is useful for two reasons. First, it gives a non-parametric form of the geometric slice condition, turning the directional curvature assumption into a fully nonlinear equation on the fiber. Second, it is the natural setting of Calabi--Bernstein type questions: one starts with a function on the fiber and asks whether the corresponding entire spacelike graph must be a slice. In this passage, obviously, one must distinguish carefully between the metric induced on the graph and the original metric of the fiber.

\vspace{1mm}

We now recall the notion of spacelike graph. Let $\Omega$ be an open subdomain of $F$, and let $u\in C^\infty(\Omega)$ satisfying $u(\Omega)\subset I$. The graph defined by $u$ in $\Mbar$ is the hypersurface
\begin{equation}\label{eq:graph}
\Sigma_u=\{(u(x),x)\, :\, x\in \Omega\}.
\end{equation}
Consider the natural embedding
\begin{equation}\label{eq:graph-embedding}
  \psi_u:\Omega\longrightarrow \Mbar,
  \qquad
  \psi_u(x)=(u(x),x),
\end{equation}
that satisfies $\psi_u(\Omega)=\Sigma_u$.
The induced metric on $\Omega$ via $\psi_u$ is
\begin{equation}\label{eq:induced-graph-metric}
  \langle \;,\; \rangle_{g_u}=-du^2+f(u)^2\langle\,,\,\rangle_F.
\end{equation}
Thus, the graph $\Sigma_u$ is spacelike if and only if
\begin{equation}\label{eq:spacelike-graph}
  |Du|_F<f(u),
\end{equation}
everywhere on $\Omega$, where $Du$ is the gradient of the function $u$ on $F$ and $|Du|_F^2:=\langle Du,Du\rangle_F$. 

\vspace{1mm}

For each spacelike graph $\Sigma_u$, we set
\begin{equation}\label{eq:W-definition}
  W=\sqrt{f(u)^2-|Du|_F^2}\,,
\end{equation}
which is smooth and positive. The unit normal vector field in the time orientation \eqref{eq:normal-orientation} is
\begin{equation}\label{eq:normal-graph}
  N=\frac{f(u)}{W}\,\Big(\partial_t+\frac{1}{f(u)^2}\,Du\Big),
\end{equation}
where $Du$ on the right-hand side means the natural lift of the gradient of $u$ to $\Mbar$. 

\vspace{1mm}

A spacelike graph $\Sigma_u$ is said to be entire if $\Omega=F$, that is,  if $u$ is defined on all of $F$. Unless explicitly stated otherwise, all spacelike graphs considered below will be assumed to be entire. In general, although $F$ is assumed geodesically complete, an entire spacelike graph may be incomplete. A well-known sufficient assumption to have completeness of a spacelike graph on $(F,\langle \, , \, \rangle_F)$ complete is, 
\begin{proposition}\label{conditions-for-completeness}
Let $\Mbar=I\times_f F$ be a GRW spacetime whose fiber $(F,g_F)$ is complete, and assume $\inf_I f>0$. If an entire spacelike graph $\Sigma_u$ satisfies
\begin{equation}\label{eq:spacelike-graph2}
  |Du|_F<\lambda f(u),  
\end{equation}
for some constant $\lambda$ with $0<\lambda<1$, then $\Sigma_u$ is complete.
\end{proposition}
\begin{proof}
From \eqref{eq:induced-graph-metric}, the Cauchy--Schwarz inequality, and \eqref{eq:spacelike-graph2}, for every $X\in T_xF$ we have
\begin{align*}
  g_u(X,X)
  &=-\langle Du,X\rangle_F^2+f(u)^2|X|_F^2\\
  &\geq\big(f(u)^2-|Du|_F^2\big)|X|_F^2\\
  &\geq(1-\lambda^2)(\inf_I f)^2|X|_F^2.
\end{align*}
Consequently, every piecewise smooth curve $\gamma$ in $F$ satisfies
\begin{equation*}
  L_{g_u}(\gamma)
  \geq\sqrt{1-\lambda^2}\,\inf_I f\,L_F(\gamma).
\end{equation*}
Since $(F,g_F)$ is complete, every divergent curve has infinite $g_F$-length and therefore infinite $g_u$-length. Hence, $(F,g_u)$ is complete.
\end{proof}
\begin{remark}\label{hyperbolic_angle} (1)
It should be noted that \eqref{eq:spacelike-graph2} means that the hyperbolic angle $\varphi$ is bounded, namely $\cosh \varphi < \frac{1}{\sqrt{1-\lambda^2}}$.  Conversely, if $\cosh \varphi < \mu$, with $\mu >1$, then we have \eqref{eq:spacelike-graph2} with $\lambda =\sqrt{1-\frac{1}{\mu^2}}$. (2)
The proof of Proposition \ref{conditions-for-completeness} shows that, for a fixed spacelike graph $\Sigma_u$, the global assumption $\inf_I f>0$ can be replaced by the weaker condition $\inf_F f(u)>0$.
\end{remark}

The next computation clarifies which metric we use. Tangent vectors to the graph are identified with tangent vectors  on $F$ through
\begin{equation}\label{eq:identification}
  X\longmapsto d\psi_u(X)\equiv X(u)\partial_t+X.
\end{equation}
For any $X,Y\in\mathfrak X(F)$, we put
\begin{equation}\label{eq:graph-B}
  \mathcal B_u(X,Y)
  =
  \Hess^Fu\,(X,Y)
  +f(u)f'(u)\langle X,Y\rangle_F
  -2\frac{f'(u)}{f(u)}\,du(X)du(Y),
\end{equation}
which is a symmetric tensor field on $F$ and satisfies
\begin{equation}\label{eq:graph-E}
  \mathcal {E}(u):= \mathcal {B}_u(Du,Du)
  =
  \Delta^F_\infty u
  +
  \frac{f'(u)}{f(u)}
  \big(f(u)^2-2|Du|_F^2\big)|Du|_F^2,
\end{equation}
The symmetric tensor field $\mathcal{B}_u$ on $F$ and (\ref{eq:identification}),
will allow writing the second fundamental form of $\Sigma_u$ on $F$ as follows.

\begin{proposition}\label{prop:graph-second-fundamental}
The Weingarten operator $A$ of a spacelike graph $\Sigma_u$ satisfies
\begin{equation}\label{eq:graph-h}
   \langle AX,Y\rangle_{g_u}
  =
  \frac{f(u)}{W}\,\mathcal B_u(X,Y),
\end{equation}
for any $X,Y\in \mathfrak{X}(F)$. Moreover, we have the following relation between the gradient of $u$ considered as a function on $\Sigma_u$ and as a function on $F$
\begin{equation}\label{eq:graph-gradient}
  \nabla^{g_u}u=\frac{1}{W^2}Du.
\end{equation}
Therefore
\begin{equation}\label{eq:graph-directional-term}
  \langle A\nabla^{g_u}u,\nabla^{g_u}u\rangle_{g_u}
  =
  \frac{f(u)}{W^5}\,\mathcal {E}(u).
\end{equation}
\end{proposition}

\begin{proof}
For each $X\in\mathfrak{X}(F)$, we write $\widetilde{X} = X(u)\partial_t + X$ for the corresponding tangent field on the graph. The connection formulas of the warped product \cite[Prop. 7.35]{ONeill1983} give
\begin{equation*}
  \nbar_{\widetilde{X}}\widetilde{Y}
  =
  \big(X(Yu)+f(u)f'(u)\langle X,Y\rangle_F\big)\partial_t
  +D_XY+\frac{f'(u)}{f(u)}\big(X(u)Y+Y(u)X\big),
\end{equation*}
where $D$ is the Levi-Civita connection of $(F,\langle\;,\; \rangle_F)$. Now, taking into account
\begin{equation*}
   \langle AX,Y\rangle_{g_u}= \langle \nbar_{\widetilde{X}}N,\widetilde{Y}\rangle
  =
  -\langle \nbar_{\widetilde{X}}\widetilde{Y},N\rangle,
\end{equation*}
where $N$ is given in \eqref{eq:normal-graph}, and we are using the Weingarten formula \eqref{eq:weingarten}, we get \eqref{eq:graph-h}.

\vspace{1mm}

For \eqref{eq:graph-gradient}, observe that for every $Y\in T_xF$,
\begin{equation*}
  \Big\langle \frac{1}{W^2}Du,Y\Big\rangle_{g_u}
  =
  \frac{1}{W^2}
  \Big(f(u)^2\langle Du,Y\rangle_F-|Du|_F^2\langle Du,Y\rangle_F\Big)
  =
  du(Y).
\end{equation*}
Combining \eqref{eq:graph-E} with \eqref{eq:graph-h} and \eqref{eq:graph-gradient}, we arrive at \eqref{eq:graph-directional-term}. This ends the proof.
\end{proof}

The following result gives the non-parametric form of the directional condition. 

\begin{corollary}\label{cor:corrected-graph-equation}
For a graph $\Sigma_u$, the conditions to be spacelike and 
\begin{equation*}
  \langle A\nabla^{g_u}u,\nabla^{g_u}u\rangle_{g_u}=0
\end{equation*}
are equivalent to
\begin{align}\tag{E}  \label{eq:graph-equation}
\begin{cases}\displaystyle
  \Delta^F_\infty u
  +\frac{f'(u)}{f(u)}\big(f(u)^2-2|Du|_F^2\big)|Du|_F^2=0, \\[4mm]
  |Du|_F<f(u).
\end{cases}
\end{align}
\end{corollary}

\subsection{Local non-rigidity of the spacelike graph equation}\label{subsec:local-nonrigidity}

The graph equation \eqref{eq:graph-equation} 
is not locally rigid by itself. The slice conclusion is a global consequence of completeness, boundedness, and the Liouville property stated in Lemma \ref{lem:liouville}. 

\begin{example}[Static case]\label{ex:static-local-nonrigidity}
Let $\Mbar=\mathbb{R}\times_c\mathbb{R}^n$ be the static GRW spacetime with $c>0$. Thus, $f\equiv c$, $f'=0$, and \eqref{eq:graph-equation} reduces to
\begin{equation*}
  \Delta^F_\infty u=0.
\end{equation*}
Every affine function
\begin{equation*}
  u(x)=\sum_{i=1}^n a_i x_i+b,
\end{equation*}
with $a_i,b\in\mathbb{R}$, satisfies $\Delta^F_{\infty} u=0$, because $\Hess^F u=0$. The corresponding graph is spacelike whenever $\sum_{i=1}^n a^2_i<c^2$. Thus, the differential equation \eqref{eq:graph-equation} alone does not force the corresponding spacelike graph to be a slice.
\end{example}

\begin{example}[One-coordinate local solutions]\label{ex:one-dimensional-local-solutions}
Let $F=\R^n$ with its Euclidean metric and $\Mbar=I\times_f \mathbb{R}^n$. Let $J\subset\mathbb{R}$ be an open interval, set $\Omega=J\times\mathbb{R}^{n-1}$, and let $v:J\to I$ be a smooth function such that $v'\neq0$ everywhere on $J$. Define $u:\Omega\to I$ by
\begin{equation*}
  u(x_1,\ldots,x_n)=v(x_1).
\end{equation*}
Then $Du=v'\frac{\partial}{\partial x_1}$, $|Du|_F^2=(v')^2$, and
$\Delta^F_\infty u=(v')^2v''$. On $\Omega$, the differential equation in \eqref{eq:graph-equation} is equivalent to
\begin{equation}\label{eq:one-dimensional-ode}
  v''+f(v)f'(v)-2\frac{f'(v)}{f(v)}(v')^2=0.
\end{equation}
Writing $y(v)=(v')^2$, so that $v''=\frac12y'(v)$ because $v'\neq0$ everywhere, we obtain
\begin{equation*}
  y'(v)-4\frac{f'(v)}{f(v)}y(v)=-2f(v)f'(v).
\end{equation*}
Multiplication by the integrating factor $f(v)^{-4}$ gives
\begin{equation*}
  \left(f(v)^{-4}y(v)\right)'=\left(f(v)^{-2}\right)'.
\end{equation*}
Hence
\begin{equation}\label{eq:one-dimensional-first-integral}
  (v')^2=f(v)^2+C f(v)^4.
\end{equation}
Choosing $C<0$ and working on a region where $1+C f(v)^2>0$, one obtains local nonconstant spacelike solutions by quadrature. More precisely, for $x_{1,0}\in J$ and $v_0=v(x_{1,0})$,
\begin{equation*}
x_1-x_{1,0}=\pm\int^v_{v_0}\frac{ds}{f(s)\sqrt{1+C f(s)^2}}.
\end{equation*}
This example again shows that the rigidity statements below are genuinely global.
\end{example}

\subsection{Global rigidity for some entire spacelike graphs}\label{subsec:global-graph-rigidity}
Now we will apply the result in Section $4$ to obtain several uniqueness results for second-order differential equations with a clear geometric meaning.  

\begin{theorem}\label{thm:graph-rigidity}
Let $f: I \rightarrow \R$ be a positive smooth function. Assume $f$ is not locally constant, has $\inf_I f>0$, and satisfies  $(\log f)'' \leq 0$. The only entire solutions $u\in C^\infty(F)$, $u(F)\subset I$, of the differential equation
\begin{align}\tag{E*} \label{eq:graph-rigidity-system}
\begin{cases}\displaystyle
            \Delta^F_\infty u
  +\frac{f'(u)}{f(u)}\big(f(u)^2-2|Du|_F^2\big)|Du|_F^2=0,\\[4mm]
            |Du|_F<\lambda f(u),
\end{cases}
\end{align}
where $\lambda$ is a constant with $0<\lambda<1$, on a complete Riemannian manifold $(F,\langle \, .\, \rangle_F)$ are the constants $u=u_0\in I$.
\end{theorem}
\begin{proof}
From the constraint in \eqref{eq:graph-rigidity-system}, the graph $\Sigma_u$ defined by an entire solution $u$ is spacelike. On the other hand, from Corollary \ref{cor:corrected-graph-equation}, the differential equation in \eqref{eq:graph-rigidity-system} is exactly the condition
\begin{equation*}
  \langle A\nabla^{g_u}u,\nabla^{g_u}u\rangle_{g_u}=0
\end{equation*}
for the height function of the graph. Finally, since $(F,g_u)$ is complete by Proposition \ref{conditions-for-completeness}, the result follows from Corollary \ref{cor:directional-zero}.
\end{proof}

\begin{remark}\label{rem:einstein-de-sitter-graphs}
For the $(n+1)$-dimensional Einstein--De Sitter spacetime in Corollary \ref{cor:einstein-de-sitter}, the warping function is $f(t)=t^{2/3}$ on $(0,\infty)$. Since $\inf_{(0,\infty)}f=0$, Theorem \ref{thm:graph-rigidity} cannot be applied directly with this whole warping interval. However, if the height satisfies $u\geq\varepsilon$ for some $\varepsilon>0$, then $u(\mathbb{R}^n)\subset I_\varepsilon=(\varepsilon/2,\infty)$ and
\begin{equation*}
 \inf_{I_\varepsilon}f=(\varepsilon/2)^{2/3}>0.
\end{equation*}
On $I_\varepsilon$, the function $f$ is strictly increasing and satisfies $(\log f)''=-2/(3t^2)<0$. Consequently, Theorem \ref{thm:graph-rigidity}, applied to this restricted interval and the complete Euclidean fiber, implies that every entire solution of \eqref{eq:graph-rigidity-system} with $f(t)=t^{2/3}$, satisfying $u\geq\varepsilon>0$ and $|Du|_{\mathbb{R}^n}<\lambda u^{2/3}$ for some constant $0<\lambda<1$, is constant. The corresponding spacelike graph is therefore a spacelike slice $\{u_0\}\times\mathbb{R}^n$, with $u_0\geq\varepsilon$.
\end{remark}

We record now the graphical form of the test function. This is the non-parametric counterpart of Proposition \ref{prop:primitive-formula}.

\begin{proposition}\label{prop:primitive-graph-formula}
Let $\Sigma_u$ be a spacelike graph in a GRW spacetime $\Mbar$ and let $P$ be defined by \eqref{eq:primitive-def}. Then, with respect to the induced metric $g_u$, one has
\begin{equation}\label{eq:primitive-graph-formula}
  \Delta^{g_u}_\infty(P(u))
  =
  \frac{f(u)^2}{W^6}
  \left(
    f(u)^3\mathcal {E}(u)
    -f'(u)|Du|_F^2W^4
  \right),
\end{equation}
where $W$ is defined in \eqref{eq:W-definition} and $\mathcal {E}(u)$ in \eqref{eq:graph-E}.
\end{proposition}

\begin{proof}
For the height function of the graph $\Sigma_u$, we have $h=u$. Hence, Proposition \ref{prop:primitive-formula} reads in this case
\begin{equation*}
  \Delta^{g_u}_\infty(P(u))
  =f(u)^2\left[f(u)\langle\alpha(\nabla^{g_u}u,\nabla^{g_u}u),\partial_t\rangle
  -f'(u)|\nabla^{g_u}u|_{g_u}^2\right].
\end{equation*}
Using \eqref{eq:normal-graph}, \eqref{eq:graph-gradient} and \eqref{eq:graph-directional-term}, we obtain
\begin{equation*}
  |\nabla^{g_u}u|_{g_u}^2=\frac{|Du|_F^2}{W^2},
  \qquad
    f(u)\langle\alpha(\nabla^{g_u}u,\nabla^{g_u}u),\partial_t\rangle
  =
  \frac{f(u)^3}{W^6}\mathcal {E}(u).
\end{equation*}
Substituting these identities into the previous formula gives \eqref{eq:primitive-graph-formula}.
\end{proof}

\begin{theorem}\label{thm:primitive-graph-consequence}
Let $f: I \to \mathbb{R}$ be a positive smooth function and assume $\inf_I f>0$.

\vspace{1mm}

\begin{enumerate}
\item[\rm (1)] If $f'\geq 0$ and $P(u)$ is bounded from below, then the only entire solutions $u\in C^\infty(F)$, $u(F)\subset I$, of the differential inequality
  \begin{align}\label{Eq2}
    \begin{cases}\displaystyle
            \Delta^F_\infty u
  +\frac{f'(u)}{f(u)}\big(f(u)^2-2|Du|_F^2\big)|Du|_F^2\leq 0,\\[4mm]
            |Du|_F<\lambda f(u),
           \end{cases}
\end{align}
where $\lambda$ is a constant with $0<\lambda<1$, on a complete Riemannian manifold $F$ are the constants $u=u_0\in I$.

\vspace*{1mm}  
  
\item[\rm (2)] If $f'\leq 0$ and $P(u)$ is bounded from above, then the only entire solutions $u\in C^\infty(F)$, $u(F)\subset I$, of the differential inequality
\begin{align}\label{Eq3}
    \begin{cases}\displaystyle
            \Delta^F_\infty u
  +\frac{f'(u)}{f(u)}\big(f(u)^2-2|Du|_F^2\big)|Du|_F^2\geq 0,\\[4mm]
            |Du|_F<\lambda f(u),
           \end{cases}
\end{align}
where $\lambda$ is a constant with $0<\lambda<1$, on a complete Riemannian manifold $F$ are the constants $u=u_0\in I$.
\end{enumerate}
\end{theorem}

\begin{proof}
(1) Since $W>0$, formula \eqref{eq:primitive-graph-formula} yields
\begin{equation*}
  \Delta^{g_u}_\infty(P(u))\leq0.
\end{equation*}
By Proposition \ref{conditions-for-completeness}, $(F,g_u)$ is complete. As $P(u)$ is bounded from below, Lemma \ref{lem:liouville} implies that $P(u)$ is constant. Since $P$ is strictly increasing, the constancy of $P(u)$ implies that $u$ is constant.

\vspace{1mm}

The case (2) is analogous. First, \eqref{eq:primitive-graph-formula} gives
\begin{equation*}
  \Delta^{g_u}_\infty(P(u))\geq0.
\end{equation*}
The boundedness from above of $P(u)$ and the Liouville property in Lemma \ref{lem:liouville} force $P(u)$ to be constant. Since $P$ is strictly increasing, it follows directly that $u$ is constant.
\end{proof}

\begin{remark}\label{rem:primitive-graph-static}
Theorem \ref{thm:primitive-graph-consequence} includes the static case $f\equiv c>0$. In this case, $P(u)=c\,(u-t_*)$ and $\mathcal{E}(u)=\Delta_\infty^F u$. Thus, on a complete Riemannian manifold $(F,\langle \; ,\, \rangle_F)$, an entire function $u\in C^\infty(F)$, with $u(F)\subset I$, is constant whenever it satisfies either
\begin{enumerate}
  \item[\rm (1)] $\Delta_\infty^F u\leq0$ and $u$ is bounded from below, or

\vspace*{1mm}
  
  \item[\rm (2)] $\Delta_\infty^F u\geq0$ and $u$ is bounded from above.
\end{enumerate}
In fact, these conclusions follow directly from Lemma \ref{lem:liouville} applied on $(F,\langle \; ,\, \rangle_F)$, so no gradient bound is needed for this analytic statement. For the geometric interpretation as a spacelike graph, one requires $|Du|_F<c$; the stronger uniform bound $|Du|_F<\lambda c$, with $0<\lambda<1$, guarantees completeness of the induced metric by Proposition \ref{conditions-for-completeness}, but is not needed for the preceding direct Liouville argument.
\end{remark}

\begin{remark}\label{rem:primitive-graph-trapped}
With the convention \eqref{eq:alpha-convention}, one has
\begin{equation*}
  \langle \alpha(\nabla^{g_u}u,\nabla^{g_u}u),\partial_t\rangle
  =\frac{f(u)^2}{W^6}\mathcal {E}(u),
\end{equation*}
from \eqref{eq:graph-directional-term}. Thus, if the spacelike graph is totally future trapped, then $\mathcal{E}(u)\leq0$; if it is totally past trapped, then $\mathcal {E}(u)\geq0$. Under the hypotheses of part {\rm(1)} of Theorem \ref{thm:primitive-graph-consequence}, the function defining a totally future trapped spacelike graph would therefore have to be constant. It would then be a spacelike slice $\{u_0\}\times F$, on which, for every nonzero tangent vector $X$,
\begin{equation*}
  \langle\alpha(X,X),\partial_t\rangle
  =\frac{f'(u_0)}{f(u_0)}|X|^2\geq 0.
\end{equation*}
This contradicts the strict negative sign required for a totally future trapped hypersurface. Consequently, there exists no entire totally future trapped spacelike graph satisfying the hypotheses of part {\rm(1)}. Analogously, under the hypotheses of part {\rm(2)}, a totally past trapped spacelike graph would have to be a slice on which
\begin{equation*}
  \langle\alpha(X,X),\partial_t\rangle
  =\frac{f'(u_0)}{f(u_0)}|X|^2\leq0,
\end{equation*}
contradicting the strict positive sign required in the totally past trapped case. Hence, there exists no entire totally past trapped spacelike graph satisfying the hypotheses of part {\rm(2)}.
\end{remark}

\begin{proposition}\label{prop:log-graph-formula}
Let $\Sigma_u$ be a spacelike graph in a GRW spacetime $\Mbar$ and set $\eta(t)=f'(t)/f(t)$. Then, we have
\begin{equation}\label{eq:log-graph-formula}
\begin{aligned}
  \Delta^{g_u}_\infty\log(f(u))
  ={}&
  \eta(u)^2\Big[
  \frac{f(u)f'(u)}{W^6}\mathcal {E}(u)
  -\eta(u)^2\frac{|Du|_F^2}{W^2}
  \\
  &\hspace{-0.5cm}
  +\big(\eta'(u)-\eta(u)^2\big)\frac{|Du|_F^4}{W^4}
  \Big],
\end{aligned}
\end{equation}
where $W$ is defined in \eqref{eq:W-definition} and $\mathcal {E}(u)$ is given in \eqref{eq:graph-E}.
\end{proposition}

\begin{proof}
For the height function of the graph, $h=u$. Applying Proposition \ref{prop:logarithmic-height}, and using formula \eqref{eq:graph-gradient} and
\begin{equation*}
  \langle\alpha(\nabla^{g_u}u,\nabla^{g_u}u),\partial_t\rangle=\frac{f(u)^2}{W^6}\,\mathcal {E}(u),
\end{equation*}
we obtain
\begin{equation*}
  \eta(u)\langle\alpha(\nabla^{g_u}u,\nabla^{g_u}u),\partial_t\rangle
  =\frac{f(u)f'(u)}{W^6}\,\mathcal {E}(u).
\end{equation*}
Substitution in \eqref{eq:log-alpha-formula} gives \eqref{eq:log-graph-formula}, ending the proof.
\end{proof}

\begin{theorem}\label{thm:log-graph-rigidity}
Let $f:I\to\mathbb{R}$ be a positive smooth function. Assume $f$ is not locally constant and $\inf_I f>0$. The only entire solutions $u\in C^\infty(F)$, $u(F)\subset I$, of the differential equation \rm{(E*)}, where $\lambda$ is a constant with $0<\lambda<1$, that satisfy $\eta'(u)\leq\eta(u)^2$ on a complete Riemannian manifold $(F,\langle \, .\, \rangle_F)$ are the constants $u=u_0\in I$.
\end{theorem}

\begin{proof}
The second condition in \rm{(E*)}
together with $\inf_I f>0$ imply, by Proposition \ref{conditions-for-completeness}, that $(F,g_u)$ is complete. The first equation in \rm{(E*)}
is equivalent to $\mathcal {E}(u)=0$. Hence, formula \eqref{eq:log-graph-formula} reduces to
\begin{equation*}
  \Delta^{g_u}_\infty\log(f(u))
  =
  \eta(u)^2\Big[
  -\eta(u)^2\frac{|Du|_F^2}{W^2}
  +\big(\eta'(u)-\eta(u)^2\big)\frac{|Du|_F^4}{W^4}
  \Big].
\end{equation*}
By the assumption on $\eta$, the right-hand side is non-positive. Since $\inf_I f>0$, the function $\log(f(u))$ is bounded from below. Applying Lemma \ref{lem:liouville} to the complete Riemannian manifold $(F,g_u)$ implies that $\log(f(u))$ is constant, and therefore $f(u)$ is constant. But $f$ is not locally constant; thus, Lemma \ref{lem:composition-constant} implies that $u$ is constant.
\end{proof}
\begin{remark}
\text{(a)} If $F$ is assumed to be compact in the previous uniqueness results, then the hyperbolic angle of the corresponding entire spacelike graph is bounded (see Remark \ref{hyperbolic_angle}). Thus, the existence of such a constant $\lambda$ is automatic. Moreover, the function $f(u)$ on $F$ has a positive global minimum value, and the function $P(u)$ is bounded from above and below. Thus, Theorems \ref{thm:graph-rigidity}, \ref{thm:primitive-graph-consequence} and \ref{thm:log-graph-rigidity} have simpler statements in this case. \text{(b)} In the differential equation \eqref{eq:graph-equation}, the second term is not signed in general, even in the case where $F$ is compact. Thus, none of the three uniqueness results can be directly deduced from Lemma \ref{lem:liouville} applied to $(F,\langle \;,\; \rangle_F)$.
\end{remark}



\begin{thebibliography}{99}
\bibitem{AlbujerAlias2009}
Albujer, A.L., Al\'{\i}as, L.J.: Spacelike hypersurfaces with constant mean curvature in the Steady State space. \emph{Proc. Amer. Math. Soc.} \textbf{137}, 711--721 (2009). \href{https://www.jstor.org/stable/20535789}{jstor.org/stable/20535789}

\bibitem{AledoRubioSalamanca2017}
Aledo, J.A., Rubio, R.M., Salamanca, J.J.: Complete spacelike hypersurfaces in generalized Robertson--Walker and the null convergence condition: Calabi--Bernstein problems. \emph{Rev. R. Acad. Cienc. Exactas Fis. Nat. Ser. A Mat. RACSAM} \textbf{111}, 115--128 (2017). \href{https://doi.org/10.1007/s13398-016-0277-3}{doi:10.1007/s13398-016-0277-3}

\bibitem{AliasCanovasColares2017}
Al\'{\i}as, L.J., C\'anovas, V.L., Colares, A.G.: Marginally trapped submanifolds in generalized Robertson--Walker spacetimes. \emph{Gen. Relativity Gravitation} \textbf{49}, 23--46 (2017). \href{https://doi.org/10.1007/s10714-017-2188-9}{doi:10.1007/s10714-017-2188-9}

\bibitem{AliasRomeroSanchez1995}
Al\'{\i}as, L.J., Romero, A., S\'anchez, M.: Uniqueness of complete spacelike hypersurfaces of constant mean curvature in generalized Robertson--Walker spacetimes. \emph{Gen. Relativity Gravitation} \textbf{27}, 71--84 (1995). \href{https://doi.org/10.1007/BF02105675}{doi:10.1007/BF02105675}

\bibitem{AliasRomeroSanchez1997}
Al\'{\i}as, L.J., Romero, A., S\'anchez, M.: Spacelike hypersurfaces of constant mean curvature and Calabi--Bernstein type problems. \emph{Tohoku Math. J.} \textbf{49}, 337--345 (1997). \href{https://doi.org/10.2748/tmj/1178225107}{doi:10.2748/tmj/1178225107}

\bibitem{AndradeCruzFigueiraLima2025}
Andrade, M., Cruz Jr., F.C., Figueira, R.F., Lima Jr., E.A.: Some rigidity results for complete totally trapped submanifolds in generalized Robertson--Walker spacetimes. \emph{Lett. Math. Phys.} \textbf{115}, 61 (2025). \href{https://doi.org/10.1007/s11005-025-01947-8}{doi:10.1007/s11005-025-01947-8}

\bibitem{AraujoMariPessoa2021}
Araujo, D.J., Mari, L., Pessoa, L.F.: Detecting the completeness of a Finsler manifold via potential theory for its infinity Laplacian. \emph{J. Differential Equations} \textbf{281}, 550--587 (2021). \href{https://doi.org/10.1016/j.jde.2021.02.005}{doi:10.1016/j.jde.2021.02.005}

\bibitem{ChoquetBruhat1976}
Choquet-Bruhat, Y.: Maximal submanifolds and submanifolds with constant mean extrinsic curvature of a Lorentzian manifold. \emph{Ann. Scuola Norm. Sup. Pisa Cl. Sci.} \textbf{3}, 361--376 (1976). \href{https://www.numdam.org/item/ASNSP-1976-4-3-3-361-0/}{numdam.org/item/ASNSP-1976-4-3-3-361-0}

\bibitem{GHV} Greub, W.H., Halperin, S., and Vanstone, R.: \emph{Connections, curvature and cohomology}. Academic Press, New York, \textbf{1}, (1972).

\bibitem{Hatcher} Hatcher, A.: \emph{Vector Bundles and K-Theory}. Version 2.2, Cornell Univ. (2017). \url{https://pi.math.cornell.edu/~hatcher/VBKT/VBpage.html} 

\bibitem{Jensen1993}
Jensen, R.: Uniqueness of Lipschitz extensions: Minimizing the sup norm of the gradient. \emph{Arch. Ration. Mech. Anal.} \textbf{123}, 51--74 (1993). \href{https://doi.org/10.1007/BF00386368}{doi:10.1007/BF00386368}

\bibitem{Juutinen2001}
Juutinen, P.: On the definition of viscosity solutions for parabolic equations. \emph{Proc. Amer. Math. Soc.} \textbf{129}, 2907--2911 (2001). \href{https://doi.org/10.1090/S0002-9939-01-05889-0}{doi:10.1090/S0002-9939-01-05889-0}

\bibitem{ONeill1983}
O'Neill, B.: \emph{Semi-Riemannian Geometry with Applications to Relativity}. Academic Press, New York (1983).

\bibitem{Penrose1965}
Penrose, R.: Gravitational collapse and space-time singularities. \emph{Phys. Rev. Lett.} \textbf{14}, 57--59 (1965). \href{https://doi.org/10.1103/PhysRevLett.14.57}{doi:10.1103/PhysRevLett.14.57}

\bibitem{RS1995} Rainer, M., Schmidt, H.J.: Inhomogeneous cosmological models with homogeneous inner hypersurface geometry. \emph{Gen. Relat. Gravit.} \textbf{27}, 1265–-1293 (1995). \href{https://doi.org/10.1007/BF02153317}{doi:10.1007/BF02153317}

\bibitem{SachsWu1977}
Sachs, R.K., Wu, H.: \emph{General Relativity for Mathematicians}. Graduate Texts in Mathematics, vol. \textbf{48}. Springer, New York (1977). \href{https://doi.org/10.1007/978-1-4612-9903-5}{doi:10.1007/978-1-4612-9903-5}

\bibitem{Senovilla2007}
Senovilla, J.M.M.: Classification of spacelike surfaces in spacetime. \emph{Class. Quantum Grav.} \textbf{24}, 3091--3124 (2007). \href{https://doi.org/10.1088/0264-9381/24/11/020}{doi:10.1088/0264-9381/24/11/020}

\bibitem{Senovilla2011} Senovilla, J.M.M.: Trapped surfaces. \emph{International Journal of Modern Physics D}, \textbf{20}, 2139--2168 (2011).
\href{https://doi.org/10.1142/S0218271811020354}{10.1142/S0218271811020354}

\end{thebibliography}
\end{document}